\documentclass[12pt,leqno,a4paper]{article}
\usepackage{amsmath,amsthm}
\usepackage[latin1]{inputenc}
\usepackage[T1]{fontenc}
\usepackage[english]{babel}
\usepackage{enumerate}

\newtheorem{theorem}{Theorem}[section]
\newtheorem{proposition}[theorem]{Proposition}

\newtheorem{corollary}[theorem]{Corollary}
\newtheorem{definition}[theorem]{Definition}

\newtheorem{remark}[theorem]{Remark}



\newcommand*\proofnamestyle{\itshape}

\begin{document}

    \title{Differential analysis of matrix convex functions II}
      \author{Frank Hansen and Jun Tomiyama}
      \date{March 9, 2007}
           
      \maketitle

      \begin{abstract}
      We continue the analysis in \cite{kn:hansen:2007:1} of matrix convex functions of a fixed order 
      defined in a real interval by differential methods as opposed to the characterization in terms
    of divided differences given by Kraus~\cite{kn:kraus:1936}. We amend and improve some points in the
    previously given presentation, and we give a number of simple but important
    consequences of matrix convexity of low orders.
    2000 AMS Classification: 26A51 and 47A63
    \end{abstract}

    \section{Introduction}

    Let $ f $ be a real function defined in an interval $ I. $ It is said to be
    $ n $-convex if
    \[
    f(\lambda A+(1-\lambda)B)\le \lambda f(A)+(1-\lambda)f(B)\qquad\lambda\in[0,1]
    \]
    for arbitrary Hermitian $ n\times n $ matrices $ A $ and $ B $ with spectra in $ I. $ It is said to be
    $ n $-concave if $ -f $ is $ n $-convex, and it is said to be $ n $-monotone if
    \[
    A\le B\quad\Rightarrow\quad f(A)\le f(B)
    \]
    for arbitrary Hermitian $ n\times n $ matrices $ A $ and $ B $ with spectra in $ I. $
    We denote by $ P_n(I) $ the set
    of $ n $-monotone functions defined in an interval $ I, $ and by $ K_n(I) $ the set of $ n $-convex functions
    defined in $ I. $
    
    We analyzed in \cite{kn:hansen:2007:1} the structure of the sets $ K_n(I) $ by differential methods and
    proved, among other things, that $ K_{n+1}(I) $ is strictly contained in $ K_n(I) $ for every natural
    number $ n. $ We discovered that some improvements of the analysis and presentation is called for, and
    this is the topic of the next section. We also noticed that the theory has quite striking applications
    for monotone or convex functions of low order, and this is exposed in the last section.
    
    \section{Improvements and amendments}
    
    \begin{definition}
    Let $ f\colon I\to\mathbf R $ be a function defined in an open interval.  
    We say that $ f $ is strictly $ n $-monotone, if $ f $ is $ n $-monotone and $ 2n-1 $ times continuously
    differentiable, and the determinant
    \[
    \det\left(\frac{f^{(i+j-1)}(t)}{(i+j-1)!}\right)_{i,j=1}^n>0
    \]
    for every $ t\in I. $ Likewise, we say that $ f $ is strictly $ n $-convex, if $ f $ is 
    $ n $-convex and $ 2n $ times continuously differentiable, and the determinant
    \[
    \det\left(\frac{f^{(i+j)}(t)}{(i+j)!}\right)_{i,j=1}^n>0
    \]
    for every $ t\in I. $    
    \end{definition}
    
    Just by inspecting the proof of \cite[Proposition 1.3]{kn:hansen:2007:1}, we realize that we 
    previously already proved the following slightly stronger result.
    
    \begin{proposition}\label{proposition: existence of strictly n-monotone and n-concave polynomials}
    Let $ I $ be a finite interval, and let $ m $ and $ n $ be natural numbers with $ m\ge 2n. $    
    There exists a strictly $ n $-concave and strictly $ n $-monotone polynomial 
    $ f_m\colon I\to\mathbf R $ of degree $ m. $
    Likewise, there exists a strictly $ n $-convex and strictly $ n $-monotone polynomial
    $ g_m\colon I\to\mathbf R $ of degree $ m. $
    \end{proposition}
    
    \begin{remark}\rm
    We would like to give some more detailed comments to the proof of \cite[Theorem 1.2]{kn:hansen:2007:1}. The theorem
    states that if $ f $ is a real $ 2n $ times continuously differentiable function defined in an open interval
    $ I, $ then the matrix
    \[
    K_n(f;t)=\begin{pmatrix}
    \displaystyle\frac{f^{i+j}(t)}{(i+j)!}
    \end{pmatrix}_{i,j=1}^n
    \]
    is positive semi-definite for each $ t\in I. $ We proved that the leading determinants of the
    matrix $ K_n(f;t) $ are non-negative for each $ t\in I. $ It is well-known that this condition is not
    sufficient to assure that the matrix itself is positive semi-definite. In the proof we wave our hands and say
    that all principal submatrices of $ K_n(f;t) $ may be obtained as a leading principal submatrix by first
    making a suitable joint permutation of the rows and columns in the Kraus matrix. But this common remedy is 
    unfortunately not working in the present situation. We therefore owe the readers to finish the proof correctly.
    \end{remark}
    
    \begin{proof}
    Let $ D_m(K_n(f;t_0)) $ for some $ t_0\in I $ denote the leading principal determinant of order $ m $ of the matrix
    $ K_n(f;t_0). $ We may according to 
    Proposition \ref{proposition: existence of strictly n-monotone and n-concave polynomials}
    choose a matrix convex function $ g $ such that 
    \[
    D_m(K_n(g;\,t_0))>0\qquad m=1,\dots,n.
    \]
    The polynomial $ p_m $ in $ \epsilon $ defined by setting
    \[
    p_m(\epsilon)=D_m(K_n(f+\epsilon g;\, t_0))
    \]
    is of degree at most $ m, $ and $ p_m(\epsilon)\ge 0 $ for
    $ \epsilon\ge 0. $ But since the coefficient to $ \epsilon^m $ in $ p_m $ is $ D_m(K_n(g;\,t_0))>0, $
    we realize that $ p_m $ is not the zero polynomial.     
    Let $ \eta_m $ be the smallest positive root of $ p_m, $ then
    \[
    p_m(\epsilon)>0\qquad 0<\epsilon<\eta_m.
    \]
    Setting $ \eta=\min\{\eta_1,\dots,\eta_n) $ we obtain
    \[
    K_n(f+\epsilon g;\, t_0)>0\qquad 0<\epsilon<\eta.
    \]
    By letting $ \epsilon $ tend to zero, we finally conclude that $ K_n(f;\, t_0) $ is positive
    semi-definite.
    \end{proof}
    
    We state in a remark after \cite[Corollary 1.5]{kn:hansen:2007:1} that the possible degrees of any polynomial
    in the gab between the matrix convex functions of order $ n $ and order $ n+1 $ defined in a finite
    interval are limited to $ 2n $ and
    $ 2n+1. $ But this is taken in the context of polynomials of degree less than or equal to $ 2n+1 $ and may be
    misunderstood. There may well be polynomials of higher degrees in the gab.
    
    \section{Scattered observations}
    
    It is well-known for which exponents the function $ t\to t^p $ is either operator monotone or operator
    convex in the positive half-axis. 
    It turns out that the same results apply if we ask for which exponents the function is
    $ 2 $-monotone or $ 2 $-convex in an open subinterval of the positive half-axis.
    
     \begin{proposition}
    Consider the function
    \[
    f(t)=t^p\qquad t\in I
    \]
    defined in any subinterval $ I $ of the positive half-axis. Then $ f $ is $ 2 $-monotone if and only if
    $ 0\le p\le 1, $ and it is $ 2 $-convex if and only if either $ 1\le p\le 2 $ or $ -1\le p\le 0. $
    \end{proposition}
    
    \begin{proof} There is nothing to prove if $ f $ is constant or linear, so we may assume that $ p\ne 0 $
    and $ p\ne 1. $
    In the first case the derivative $ f'(t)=p t^{p-1} $ should be non-negative so $ p> 0, $ 
    and it may be written \cite[Chapter VII Theorem IV]{kn:donoghue:1974} on the form
    \[
    f'(t)=\frac{1}{c(t)^2}\qquad t\in I
    \]
    for $ c(t)=p^{-1/2} t^{(1-p)/2} $ and this function is concave only for $ 0< p\le 1. $
    One may alternatively consider the determinant
    \[
    \begin{array}{rl}
    \det\begin{pmatrix}
        f'(t) & \displaystyle\frac{f''(t)}{2!}\\[2ex]
        \displaystyle\frac{f''(t)}{2!} & \displaystyle\frac{f^{(3)}(t)}{3!}
        \end{pmatrix}
    &=\det\begin{pmatrix}
        p t^{p-1} & \displaystyle\frac{p(p-1)t^{p-2}}{2}\\[2ex]
        \displaystyle\frac{p(p-1)t^{p-2}}{2} 
        & \displaystyle\frac{p(p-1)(p-2) t^{p-3}}{6}
        \end{pmatrix}\\[2ex]
    &=\displaystyle -\frac{1}{12}p^2 (p-1)(p+1) t^{2p-4}
    \end{array}
    \] 
    and note that the matrix is positive semi-definite only for $ 0\le p\le 1. $
    
    The second derivative may be written \cite[Theorem 2.3]{kn:hansen:2007:1} on the form
    \[
    f''(t)=p(p-1) t^{p-2}=\frac{1}{d(t)^3}\qquad t\in I
    \]
    for $ d(t)=(p(p-1))^{-1/3} t^{(2-p)/3}, $ and this function is concave only for $ -1\le p< 0 $
    or $ 1< p\le 2. $ One may alternatively consider the determinant
    \[
    \begin{array}{l}
    \det\begin{pmatrix}
        \displaystyle\frac{f''(t)}{2!} & \displaystyle\frac{f^{(3)}(t)}{6}\\[2ex]
        \displaystyle\frac{f^{(3)}(t)}{6} & \displaystyle\frac{f^{(4)}(t)}{24}
        \end{pmatrix}\\[6ex]
    =\det\begin{pmatrix}
        \displaystyle\frac{p(p-1) t^{p-2}}{2} & \displaystyle\frac{p(p-1)(p-2)t^{p-3}}{6}\\[2ex]
        \displaystyle\frac{p(p-1)(p-2)t^{p-3}}{6} 
        & \displaystyle\frac{p(p-1)(p-2)(p-3) t^{p-4}}{24}
        \end{pmatrix}\\[6ex]
    =\displaystyle -\frac{1}{144}p^2 (p-1)^2 (p-2)(p+1) t^{2p-6}
    \end{array}
    \] 
    and note that the matrix is positive semi-definite only for $ -1\le p\le 0 $ or $ 1\le p\le 2. $
    \end{proof}
    
    The observation that the function $ t\to t^p $ is $ 2 $-monotone only for $ 0\le p\le 1 $ has appeared
    in the literature in different forms, cf. \cite[1.3.9 Proposition]{kn:pedersen:1979} or
    \cite{kn:ji:2003}.
    
    It is known that the derivative of an operator monotone function defined on an infinite interval 
    $ (\alpha,\infty) $ is completely monotone \cite[Page 86]{kn:donoghue:1974}. We give a parallel
    result for matrix monotone functions implying this observation, and
    extend the analysis to matrix convex functions. 
    
    \begin{theorem}
    Consider a function $ f $ defined in an interval of the form $ (\alpha,\infty) $ for some real $ \alpha. $
    \begin{enumerate}[(i)]
    
    \item If $ f $ is $ n $-monotone and $ 2n-1 $ times continuously differentiable, then
    \[
    (-1)^k f^{(k+1)}(t)\ge 0\qquad k=0,1,\dots,2n-2.
    \]
    Therefore, the function $ f $ and its even derivatives up to order $ 2n-4 $ are concave functions, and 
    the odd derivatives up to order $ 2n-3 $ are convex functions.
    
    \item If $ f $ is $ n $-convex and $ 2n $ times continuously differentiable, then
    \[
    (-1)^k f^{(k+2)}(t)\ge 0\qquad k=0,1,\dots,2n-2.
    \]
    Therefore, the function $ f $ and its even derivatives up to order $ 2n-2 $ are convex functions, and 
    the odd derivatives up to order $ 2n-3 $ are concave functions.
    \end{enumerate} 
    
    \end{theorem}
    
    \begin{proof}
    We may assume $ n\ge 2. $ To prove the first assertion we may 
    write \cite[Chapter VII Theorem IV]{kn:donoghue:1974} the derivative $ f' $ on the form
    \[
    f'(t)=\frac{1}{c(t)^2},
    \]
    where $ c $ is a positive concave function. Since $ c $ is defined in an infinite interval
    it has to be increasing, therefore $ f' $ is decreasing and thus $ f''\le 0. $
    Since $ f $ is $ n $-monotone, it follows from Dobsch' condition \cite{kn:Dobsch:1937} 
    that the odd derivatives satisfy 
    \[
    f^{(2k+1)}\ge 0\qquad k=0,1,\dots, n-1.
    \]
    The odd derivatives $ f^{(2k+1)} $ are thus convex for $ k=0,1,\dots,n-2. $ If the third derivative
    $ f^{(3)}, $ which is a convex function, were strictly increasing at any point, then it would go towards
    infinity and the second derivative would eventually be positive for large $ t. $ But this contradicts
    $ f''\le 0, $ so $ f^{(3)} $ is decreasing and thus the fourth derivative $ f^{(4)}\le 0. $ This argument
    may now be continued to prove the first assertion.
    
    To prove the second assertion we may write \cite[Theorem 2.3]{kn:hansen:2007:1} the second derivative $ f'' $
    on the form
    \[
    f''(t)=\frac{1}{d(t)^3},
    \]
    where $ d $ is a positive concave function. Since $ d $ is defined in an infinite interval
    it has to be increasing, therefore $ f'' $ is decreasing and thus $ f^{(3)}\le 0. $
    Since $ f $ is $ n $-convex, it follows \cite[Theorem 1.2]{kn:hansen:2007:1} that the even derivatives satisfy 
    \[
    f^{(2k)}\ge 0\qquad k=1,\dots, n.
    \]
    The statement now follows in a similar way as for the first assertion.    
    \end{proof}
    
    \begin{corollary}
    The second derivative of an operator convex function defined in an infinite interval $ (\alpha,\infty) $
    is completely monotone.
    \end{corollary}

    \begin{remark}
    The indefinite integral $ g(t)=\int f(t)\,dt $ of a $ 2 $-monotone function $ f $ is $ 2 $-convex.
    \end{remark}
    
    \begin{proof} The second derivative may be written on the form
    \[
    g''(t)=f'(t)=\frac{1}{c(t)^2}=\frac{1}{(c(t)^{2/3})^3}
    \]
    for some positive concave function $ c. $ Since the function $ t\to t^{2/3} $ is increasing and concave, we
    conclude that $ t\to c(t)^{2/3} $ is concave. The statement then follows from the characterization
    of $ 2 $-convexity.
    \end{proof}
    
    It is known in the literature that operator monotone or operator convex functions defined in the whole
    real line are either affine or quadratic, and this fact is established by appealing to the 
    representation theorem of Pick functions. But the situation is far more general, and the results only
    depend on monotonicity or convexity on two by two matrices.
    
    \begin{theorem}
    Let $ f $ be a function defined in the whole real line. If $ f $ is $ 2 $-monotone and three times
    continuously differentiable, then it is necessarily affine. If $ f $ is $ 2 $-convex and four times
    continuously differentiable, then it is necessarily quadratic.    
    \end{theorem}    
    
    \begin{proof} In the first case the derivative $ f' $ may
    be written \cite[Chapter VII Theorem IV]{kn:donoghue:1974} on the form $ f'(t)=c(t)^{-2} $ for
    some positive concave function $ c $ defined in the real line, while in the second case 
    the second derivative $ f'' $ may be written
    \cite[Theorem 2.3]{kn:hansen:2007:1}  on the form $ f''(t)=d(t)^{-3} $ for some positive concave function
    $ d $ defined in the real line. The assertions now follows since 
    a positive concave function defined in the whole real line necessarily is constant.    
    \end{proof}

{\footnotesize

   

      \vfill

      \noindent Frank Hansen: Department of Economics, University
       of Copenhagen, Stu\-die\-straede 6, DK-1455 Copenhagen K, Denmark.\\[1ex]

       \noindent Jun Tomiyama: Department of Mathematics and Physics. Japan Women's University.
       Mejirodai Bunkyo-ku, Tokyo, Japan.

       }

      \end{document}